\documentclass[11pt, twoside]{article}
\usepackage{amsmath,amssymb}
\usepackage{float}
\usepackage{multirow}
\usepackage{graphicx}
\graphicspath{{fig_eps/}}
\usepackage{epsfig}
\usepackage{cite}
\usepackage{cases}
\usepackage{epstopdf}
\usepackage{caption}
\usepackage{subfigure}
\usepackage{titlesec}
\usepackage{textcomp}
\usepackage{textcomp,booktabs}
\usepackage{authblk}
\usepackage[usenames,dvipsnames]{color}
\usepackage{colortbl}
\definecolor{mygray}{gray}{0.9}
\definecolor{mypink}{rgb}{0.99,0.91,0.95}
\definecolor{mycyan}{cmyk}{0.3,0,0,0}
\usepackage{booktabs}
\usepackage[margin=2cm]{geometry}
\usepackage[colorlinks,linkcolor=red,anchorcolor=green,citecolor=blue,CJKbookmarks=True]{hyperref}
\newtheorem{Theorem}{Theorem}[section]

\allowdisplaybreaks
\begin{document}

\title{Backward bifurcation arising from decline of immunity against emerging infectious diseases}
\author[1]{Shuanglin Jing\thanks{shuanglinjing@hrbeu.edu.cn}}
\author[2]{Ling Xue\thanks{lxue@hrbeu.edu.cn}}
\author[2]{Jichen Yang\thanks{jichen.yang@hrbeu.edu.cn}}
\affil[1]{Department of Mathematics, Lanzhou Jiaotong University, Lanzhou 730070, China}
\affil[2]{College of Mathematical Sciences, Harbin Engineering University, Harbin 150001, China}

\date{\today}
\maketitle
\begin{abstract}
Decline of immunity is a phenomenon characterized by immunocompromised host and plays a crucial role in the epidemiology of emerging infectious diseases (EIDs) such as COVID-19. In this paper, we propose an age-structured model with vaccination and reinfection of immune individuals. We prove that the disease-free equilibrium of the model undergoes backward and forward transcritical bifurcations at the critical value of the basic reproduction number for different values of parameters. We illustrate the results by numerical computations, and also find that the endemic equilibrium exhibits a saddle-node bifurcation on the extended branch of the forward transcritical bifurcation. These results allow us to understand the interplay between the decline of immunity and EIDs, and are able to provide strategies for mitigating the impact of EIDs on global health.
\end{abstract}
\section{Introduction}
The emergence of emerging infectious diseases (EIDs), such as COVID-19, has underscored the critical need for a deeper understanding of global health challenges. These diseases, with their rapid spread and potential for high mortality, have exposed vulnerabilities in our current health systems and highlighted the importance of studying the immune response to pathogens \cite{xue2022infectivity,jing2023vaccine}. To effectively control these diseases, we need to thoroughly understand the immune response against the pathogens involved. In ongoing research, a phenomenon known as backward bifurcation has attracted more and more attention, particularly in the context of decline of immunity associated with EIDs \cite{wangari2024emergence}. In mathematical epidemiology, backward bifurcation occurs when the basic reproduction number ($\mathcal{R}_0$) falls below 1, revealing a more complex dynamics where a small unstable endemic equilibrium emerges alongside locally asymptotically stable disease-free and larger endemic equilibria. This complexity highlights the nonlinear nature of disease transmission and the intricate interplay between immunological factors and disease dynamics \cite{martcheva2015introduction,martcheva2020lyapunov,castillo2004dynamical,yang2022backward}. The decline of immunity against EIDs is a multifaceted process that can stem from the natural waning of immune responses over time following infection or the evolution of pathogens to evade immune detection \cite{xue2022infectivity,jing2023vaccine}. This decline in immunity poses a significant threat, as it can lead to reinfection of previously immune individuals, particularly those with compromised immune systems \cite{rahman2022covid}. Understanding the mechanisms behind this decline and its implications for disease transmission is crucial for developing effective control strategies.

The theoretical framework of backward bifurcation offers a powerful tool for studying these dynamics. By modeling the effects of changes in immunological parameters on disease transmission, we can gain insights into how different factors contribute to the bifurcation of equilibrium states. This understanding can then be used to predict disease trends, evaluate the impact of interventions, and optimize resource allocation. The significance of this research extends beyond the immediate need to control current outbreaks. It provides a foundation for developing more effective and sustainable strategies for preventing and controlling future EIDs. By enhancing our understanding of the immune response and disease transmission dynamics, we can build more resilient health systems that are better prepared to face the challenges posed by emerging infectious diseases.

In this paper, we establish an age-structured model, which takes into account the reinfection of immune individuals due to the decline of immunity, and has the form
\begin{equation}\label{COVID_BB_EQ1}
\left\{
\begin{split}
\frac{{\rm d} S}{{\rm d} t}\ &=\ \Lambda-\beta_{s}S\frac{\theta E+\epsilon A+I}{N} -(\alpha+u)S,\\
\frac{{\rm d} E}{{\rm d} t}\ &=\ \left(\beta_{s}S+\int^{\infty}_{0}\beta_{r}(\tau)r(t,\tau){\rm d}\tau\right)\frac{\theta E+\epsilon A+I}{N}-(\sigma+u)E,\\
\frac{{\rm d} A}{{\rm d} t}\ &=\ (1-\rho)\sigma E-(\gamma_{A}+u)A,\\
\frac{{\rm d} I}{{\rm d} t}\ &=\ \rho\sigma E-(\gamma_{I}+\mu+u)I,\\
\frac{\partial r}{\partial t}+\frac{\partial r}{\partial \tau}\ &=\ -\beta_{r}(\tau)\frac{\theta E+\epsilon A+I}{N} r-ur,
\end{split}
\right.
\end{equation}
with the boundary condition
\begin{equation}
r(t,0)\ =\ \alpha S+\gamma_{A}A+\gamma_{I}I,\label{COVID_BB_EQ1_1}
\end{equation}
where the real-valued functions $S=S(t)$, $E=E(t)$, $A=A(t)$, and $I=I(t)$ denote the number of susceptible, latent, asymptomatic infected, and symptomatic infected individuals at time $t\geq0$, respectively, and $r=r(t,\tau)$ denotes the density of immune individuals with immune age $\tau\geq0$ at time $t$. The total population at time $t$ is denoted by $N(t):=S(t)+E(t)+A(t)+I(t)+\int^{\infty}_{0}r(t,\tau)\mbox{d}\tau$. We denote by $\Lambda>0$ the recruitment rate of susceptible individuals; $\mu,u>0$ are the death rate of symptomatic individuals and the natural death rate, respectively; $\sigma>0$ is the rate at which latent individuals progress to the next stage; $\gamma_{A},\gamma_{I}>0$ are the recovery rates of asymptomatic and symptomatic individuals, respectively; $\rho>0$ quantifies the proportion of symptomatic infected individuals; $\alpha>0$ is the vaccination rate of susceptible individuals; $\theta, \epsilon\in(0,1]$ are the coefficients for reduced transmission probabilities of latent and asymptomatic infected individuals, respectively; $\beta_{s}>0$ is the transmission rate of latent and infected individuals infecting susceptible individuals. Due to the gradual decline of immunity \cite{risk2022covid}, the immune individuals are prone to be reinfected by the latent and infected individuals with the rate $\beta_{r}(\tau)(\theta E+\epsilon A+I)/N$ at the stage $\tau$. Hence, the dynamics of immune individuals can be described by the hyperbolic PDE (i.e., the last equation) in \eqref{COVID_BB_EQ1}, where $\beta_{r}(\tau)>0$ is the transmission rate, and it is monotonically increasing and bounded, and thus $\beta_{r}\in L^{\infty}(0,\infty)$.

It is easy to prove that the solutions of \eqref{COVID_BB_EQ1} are non-negative with non-negative initial data. With the boundary condition \eqref{COVID_BB_EQ1_1}, the disease-free equilibrium, denoted by $P_{0}$, is given by
\[
P_0:=\left(S^{0},E^{0},A^{0},I^{0},r^{0}(\tau)\right) = \left(\frac{\Lambda}{\alpha+u},\ 0,\ 0,\ 0,\ \frac{\alpha\Lambda}{\alpha+u}{\rm e}^{-u\tau}\right).
\]
The basic reproduction number (cf. \cite{van2002reproduction,diekmann1990definition} for the details of computations) $\mathcal{R}_{0}$ is given by
\begin{equation*}
\mathcal{R}_{0}=\frac{1}{(\sigma+u)N^{0}}\left(\beta_{s}S^{0}+\int^{\infty}_{0}\beta_{r}(\tau)r^{0}(\tau){\rm d}\tau\right)\left(\theta+\frac{(1-\rho)\sigma\epsilon}{\gamma_{A}+u}+\frac{\rho\sigma}{\gamma_{I}+\mu+u}\right),
%\label{EQ_DI3}
\end{equation*}
where $N^{0}:=S^{0}+\int^{\infty}_{0}r^{0}(\tau){\rm d}\tau$. The rescaling $\beta_{r}(\tau)=\bar{\beta}\beta_{r0}(\tau)$ allows us to select the bifurcation parameter $\bar{\beta}>0$ with $\beta_{r0}(\tau)\in(0,1]$, and we set $\bar{\beta}^{*}$ is the value of $\bar{\beta}$ such that $\mathcal{R}_{0}= 1$.

In this paper, we study the bifurcations of \eqref{COVID_BB_EQ1}. Using Lyapunov-Schmidt reduction (cf., e.g., \cite{kielhofer2012bifurcation,martcheva2020lyapunov}), we analytically prove the backward and forward transcritical bifurcations from $P_{0}$ at $\mathcal{R}_{0}=1$ for different values of parameters. In preparation for the main result of this paper, we define the following quantity which characterizes the quadratic nonlinearity of \eqref{COVID_BB_EQ1},
\begin{align}\label{EQ_DI4}
a :=\ & \frac{2u}{\Lambda}\left(\theta+\frac{(1-\rho)\sigma\epsilon}{\gamma_{A}+u}+\frac{\rho\sigma}{\gamma_{I}+\mu+u}\right) \Bigg[\bar{\beta}^{*}\left(\frac{\alpha(\sigma+u)}{\alpha+u}+\frac{\gamma_{A}(1-\rho)\sigma}{\gamma_{A}+u}+\frac{\gamma_{I}\rho\sigma}{\gamma_{I}+\mu+u}\right)\int^{\infty}_{0}\beta_{r0}(\tau){\rm e}^{-u\tau}{\rm d}\tau\nonumber\\
&\quad -\frac{u}{\alpha+u}\left(\beta_{s}+\alpha\bar{\beta}^{*}\int^{\infty}_{0}\beta_{r0}(\tau){\rm e}^{-u\tau}{\rm d}\tau\right)\left(1+\frac{(1-\rho)\sigma}{u}+\frac{(\gamma_{I}+u)\rho\sigma}{\gamma_{I}+\mu+u}\right)\\
&\quad -\frac{\alpha u}{\alpha+u}(\bar{\beta}^{*})^{2}\left(\theta+\frac{(1-\rho)\sigma\epsilon}{\gamma_{A}+u}+\frac{\rho\sigma}{\gamma_{I}+\mu+u}\right)\int^{\infty}_{0}\beta_{r0}(\tau){\rm e}^{-u\tau}\int^{\tau}_{0}\beta_{r0}(h)\,{\rm d}h\,{\rm d}\tau\Bigg],\nonumber%\\
\end{align}
and our main theorem is as follows.
\begin{Theorem}\label{t:bif}
At $\mathcal{R}_{0}=1$, model \eqref{COVID_BB_EQ1} exhibits a backward transcritical bifurcation for $a>0$, and a forward transcritical bifurcation for $a<0$.
\end{Theorem}

We remark that the quantity $a$ defined as in  \eqref{EQ_DI4} is continuous in the parameters, thus $a$ can be zero. In such case, the transcritical bifurcation may degenerate into pitchfork bifurcations. However, this is beyond the scope of this paper and we do not pursue this further.

We illustrate the results in Theorem~\ref{t:bif} with some numerical computations, which also suggest that the forward bifurcating branch extends to the point $\mathcal{R}_{0} = \mathcal{R}_0^{\rm bsn}>1$, at which the endemic equilibrium undergoes a saddle-node bifurcation; moreover, the bistable state, i.e., the coexistence of two stable endemic equilibria, occurs for some values of parameters.

This paper is organized as follows. In section~\ref{sec2}, we give a proof of Theorem~\ref{t:bif}. We present some numerical computations in section~\ref{sec3} and provide a short discussion in section~\ref{sec4}.

\section{Bifurcation analysis}\label{sec2}
We consider the stationary solutions to \eqref{COVID_BB_EQ1} with the boundary condition \eqref{COVID_BB_EQ1_1}. The proof is essentially based upon Lyapunov-Schmidt reduction, we refer the readers to, e.g., \cite{kielhofer2012bifurcation,martcheva2020lyapunov}, for more details.

We define the function space $X:=\mathbb{R}^{4}\times L^{1}(0,\infty)$, and the nonlinear operator $\mathcal{F}: X\times(0,\infty)\rightarrow X$ as follows
\begin{equation*}
\mathcal{F}(\phi,\bar{\beta})=
\begin{pmatrix}
\Lambda-\beta_{s}S(\theta E+\epsilon A+I)/N -(\alpha+u)S \\
(\beta_{s}S+\bar{\beta}\int^{\infty}_{0}\beta_{r0}(\tau)r(t,\tau){\rm d}\tau)(\theta E+\epsilon A+I)/N-(\sigma+u)E \\
(1-\rho)\sigma E-(\gamma_{A}+u)A\\
\rho\sigma E-(\gamma_{I}+\mu+u)I\\
-\partial_{\tau}r-\bar{\beta}\beta_{r0}(\tau)r(\theta E+\epsilon A+I)/N-ur
\end{pmatrix}
\end{equation*}
with
$\phi=(S,E,A,I,r)^{\mathrm{T}}\in X$ and $\bar\beta\in(0,\infty)$. Linearizing $\mathcal{F}(\phi,\bar{\beta})$ in the disease-free equilibrium $\phi=P_0$ and evaluating at the parameter value $\bar\beta=\bar\beta^*$ gives the linear operator $\mathcal{A}:=\mathcal{F}_{\phi}(P_{0},\bar{\beta}^{*}):\mathcal{D}(\mathcal{A})\subset X\rightarrow X$ with the domain
\[
\mathcal{D}(\mathcal{A}):=\mathbb{R}^{4}\times W^{1,1}(0,\infty),
\]
which has the form
\begin{equation*}
\mathcal A =
\begin{pmatrix}
-(\alpha+u) & -\theta D & -\epsilon D & -D & 0 \\
0 & \theta B - (\sigma+u) & \epsilon B & B & 0 \\
0 & (1-\rho)\sigma & -(\gamma_A + u) & 0 & 0 \\
0 & \rho\sigma & 0 & -(\gamma_I+\mu+u) & 0 \\
0 & -\theta \mathcal E(\tau) & -\epsilon\mathcal E(\tau) & -\mathcal E(\tau) & -(\partial_\tau+u)
\end{pmatrix},
\end{equation*}
where $B:=(\beta_{s}S^{0}+\bar{\beta}^{*}\int^{\infty}_{0}\beta_{r0}(\tau)r^{0}(\tau){\rm d}\tau)/N^{0}$, $D:=\beta_s S^0/N^0$, $\mathcal E(\tau):=\bar\beta^*\beta_{r0}(\tau)r^0(\tau)/N^0$. It is easy to verify that $\mathcal A$ has a simple zero eigenvalue at $\mathcal R_0=1$ by considering the characteristic equation derived from the eigenvalue problem $\mathcal A x=\lambda x$ for non-trivial $x=(x_1,x_2,x_3,x_4,x_5)^{\mathrm T}\in \mathcal D(\mathcal A)$, subject to the boundary condition
\begin{equation}\label{e:bnd}
x_{5}(0)=\alpha x_{1}+\gamma_{A}x_{3}+\gamma_{I}x_{4},
\end{equation}
which stems from the linearization of \eqref{COVID_BB_EQ1_1} in $P_0$. Solving the equations $\mathcal A x=0$ with \eqref{e:bnd} gives the basis, denoted by $\hat x= (\hat x_{1},\hat x_{2},\hat x_{3},\hat x_{4},\hat x_{5})^{\mathrm{T}}\in \mathcal D(\mathcal A)$, of the kernel of $\mathcal A$, i.e., $\ker(\mathcal A) = {\rm span}\{\hat x\}$, whose elements take the form
%Hence, w}e can obtain the eigenfunction of $\mathcal A$ with the corresponding zero eigenvalue by solving the equations $\mathcal A x=0$ with \eqref{e:bnd}. {\JYc We denote by $\hat x= (\hat x_{1},\hat x_{2},\hat x_{3},\hat x_{4},\hat x_{5})^{\mathrm{T}}\in\mathcal D(\mathcal A)$ the basis of} the kernel of $\mathcal A$, %is spanned by $\hat x$,
%i.e., $\ker(\mathcal A) = {\rm span}\{\hat x\}$, {\JYc whose elements take} the form %$\hat{x} = (\hat x_{1},\hat x_{2},\hat x_{3},\hat x_{4},\hat x_{5})^{\mathrm{T}}$ with
\begin{align*}
\hat x_{1}&=-\frac{\beta_{s}S^{0}}{(\alpha+u)N^{0}}\left(\theta+\frac{(1-\rho)\sigma\epsilon}{\gamma_{A}+u}+\frac{\rho\sigma}{\gamma_{I}+\mu+u}\right),\quad \hat x_2=1,\quad \hat x_{3}=\frac{(1-\rho)\sigma}{\gamma_{A}+u}, \quad \hat x_{4}=\frac{\rho\sigma}{\gamma_{I}+\mu+u},\\
\hat x_5&=\hat x_{5}(\tau)=\left(\hat x_{5}(0)-\frac{\bar\beta^*}{N^0}\left(\theta+\frac{(1-\rho)\sigma\epsilon}{\gamma_{A}+u}+\frac{\rho\sigma}{\gamma_{I}+\mu+u}\right)
\int^{\tau}_{0}\beta_{r0}(h)r^{0}(h){\rm e}^{uh}{\rm d}h\right){\rm e}^{-u\tau},
\end{align*}
where the boundary value $\hat x_5(0)= \alpha\hat x_1 + \gamma_A\hat x_3 + \gamma_I \hat x_4$.

Next, we discuss the adjoint operator of $\mathcal A$, denoted by $\mathcal A^*$. Based upon the Riesz representation theorem on the identification of the dual space of $L^1$, we choose the function space $X^*:=\mathbb{R}^{4}\times L^{\infty}(0,\infty)$ and the domain of $\mathcal{A}^{*}$ as follows
\[
\mathcal{D}(\mathcal{A}^{*}):=\mathbb{R}^{4}\times W^{1,\infty}(0,\infty)\subset X^*.
%\mathcal{D}(\mathcal{A}^{*}):=\big\{\xi=(\xi_1,\xi_2,\xi_3,\xi_4,\xi_5)\in\mathbb{R}^{4}\times W^{1,\infty}(0,\infty): \xi_5(\infty)=0\big\}\subset X^*.
\]
It is well-known that for all $x\in\mathcal D(\mathcal A)$ and $\xi\in\mathcal D(\mathcal A^*)$, the adjoint operator $\mathcal A^*$ is unique and satisfies $\langle\mathcal{A}x,\xi\rangle=\langle x,\mathcal{A}^{*}\xi\rangle$, where the bilinear form $\langle h,g \rangle:=\sum_{j=1}^{4}h_j g_j + \int_0^\infty h_5 g_5\,{\rm d}\tau$ for any $h=(h_1,\dots,h_5)\in X$ and $g=(g_1,\dots,g_5)\in X^*$.
Hence, for all $\xi=(\xi_1,\xi_2,\xi_3,\xi_4,\xi_5)^{\mathrm T}\in\mathcal D(\mathcal A^*)$ we have
% using integration by parts with the {\JYc boundary} condition $x_{5}(0)=\alpha x_{1}+\gamma_{A}x_{3}+\gamma_{I}x_{4}$, and $\xi_5(\infty)=0$, and {\JYc the fact that} $\lim_{\tau\to\infty}x_5(\tau)=0$ for $x_5\in W^{1,1}(0,\infty)$,
\begin{align*}
\langle\mathcal{A}x,\xi\rangle\ =\ & -\left[D(\theta x_{2}+\epsilon x_{3}+x_{4})+(\alpha+u)x_{1}\right]\xi_{1} +\left[B(\theta x_{2}+\epsilon x_{3}+x_{4})-(\sigma+u)x_{2}\right]\xi_{2}\\
&\ +\left[(1-\rho)\sigma x_{2}-(\gamma_{A}+u)x_{3}\right]\xi_{3}+\left[\rho\sigma x_{2}-(\gamma_{I}+\mu+u)x_{4}\right]\xi_{4}\\
&\ -\int^{\infty}_{0}\left[\mathcal E(\tau)(\theta x_{2}+\epsilon x_{3}+x_{4})+(\partial_{\tau}+u)x_{5}(\tau)\right]\xi_{5}(\tau)\,{\rm d}\tau\\
=\ & -(\alpha+u)\xi_{1}x_{1} +\left[-\theta D\xi_{1}+\theta B\xi_{2}-(\sigma+u)\xi_{2}+(1-\rho)\sigma\xi_{3}+\rho\sigma\xi_{4}-\theta\mathcal{C}(\xi_5)\right]x_{2}\\
&\ +\left[-\epsilon D\xi_{1} + \epsilon B\xi_{2}-(\gamma_{A}+u)\xi_{3}-\epsilon\mathcal{C}(\xi_5)\right]x_{3} +\left[-D\xi_{1}+ B\xi_{2} -(\gamma_{I}+\mu+u)\xi_{4}- \mathcal{C}(\xi_5)\right]x_{4}\\
&\ - x_5(\infty)\xi_5(\infty)+ x_5(0)\xi_5(0) +\int^{\infty}_{0}\big[\partial_{\tau}\xi_{5}(\tau)-u\xi_{5}(\tau)\big]x_{5}(\tau)\,{\rm d}\tau\\
%=\ & \left[-(\alpha+u)\xi_{1}+\alpha\xi_5(0)\right]x_{1} +\left[-\theta D\xi_{1}+\theta B\xi_{2}-(\sigma+u)\xi_{2}+(1-\rho)\sigma\xi_{3}+\rho\sigma\xi_{4}-\theta\mathcal{C}(\xi_5)\right]x_{2}\\
%&\ +\left[-\epsilon D\xi_{1} + \epsilon B\xi_{2}-(\gamma_{A}+u)\xi_{3}-\epsilon\mathcal{C}(\xi_5)+\gamma_A\xi_5(0)\right]x_{3}\\
%&\ +\left[-D\xi_{1}+ B\xi_{2} -(\gamma_{I}+\mu+u)\xi_{4}- \mathcal{C}(\xi_5)+\gamma_I\xi_5(0)\right]x_{4}\\
%&\ +\int^{\infty}_{0}\big[\partial_{\tau}\xi_{5}(\tau)-u\xi_{5}(\tau)\big]x_{5}(\tau)\,{\rm d}\tau\\
=\ &\langle x, \mathcal A^*\xi \rangle, %-x_{5}(\infty)\xi_{5}(\infty),
\end{align*}
where the linear functional
$\mathcal{C}(\xi_5):=\int_0^\infty \mathcal E(\tau)\xi_5(\tau)\,{\rm d}\tau$. Combining the condition \eqref{e:bnd} and the fact that $\lim_{\tau\to\infty}x_5(\tau)=0$ for $x_5\in W^{1,1}(0,\infty)$, the adjoint operator $\mathcal{A}^{*}:\mathcal{D}(\mathcal{A}^{*})\subset X^{*}\rightarrow X^{*}$ takes the form
\begin{equation*}
\mathcal A^* =
\begin{pmatrix}
-(\alpha+u) & 0 & 0 & 0 & 0 \\
-\theta D & \theta B-(\sigma+u) & (1-\rho)\sigma & \rho\sigma & -\theta \mathcal C(\cdot) \\
-\epsilon D & \epsilon B & -(\gamma_A+u) & 0 & -\epsilon\mathcal C(\cdot) \\
-D & B & 0 & -(\gamma_I+\mu+u) & -\mathcal C(\cdot) \\
0 & 0 & 0 & 0 & \partial_\tau-u
\end{pmatrix}
\end{equation*}
subject to the adjoint boundary conditions
\[
\xi_5(0) = \xi_5(\infty) = 0.
\]
%\begin{equation*}
%\mathcal{A}^{*}\xi=
%\begin{pmatrix}
%-(\alpha+u)\xi_{1}+\alpha \xi_{5}(0)\\
%-\theta D\xi_{1}+ \theta B\xi_{2}-(\sigma+u)\xi_{2}+(1-\rho)\sigma\xi_{3}+\rho\sigma\xi_{4}
%-\theta\mathcal{C}(\xi_{5})\\
%-\epsilon D\xi_{1}+\epsilon B\xi_{2}-(\gamma_{A}+u)\xi_{3}-\epsilon\mathcal{C}(\xi_{5})+\gamma_A \xi_{5}(0)\\
%-D\xi_{1}+B\xi_{2}-(\gamma_{I}+\mu+u)\xi_{4}-\mathcal{C}(\xi_{5})+\gamma_I \xi_{5}(0)\\
%\partial_{\tau}\xi_{5}-u\xi_{5}
%\end{pmatrix}.
%\end{equation*}
Solving the equations $\mathcal A^*\xi=0$ with such boundary conditions
%{\JYc adjoint boundary} condition $\xi_5(0)=0$ ({\JYc note that} $\xi_5\notin L^\infty(0,\infty)$ for any non-zero {\JYc boundary values})
gives the eigenfunction in the kernel of $\mathcal A^*$, denoted by $\hat\xi$, which has the form
\[
\hat{\xi}=\left(0,\ 1,\ \frac{\epsilon B}{\gamma_{A}+u},\ \frac{B}{\gamma_{I}+\mu+u},\ 0\right)^{\mathrm{T}}.
\]

Differentiating $\mathcal F$ twice with respect to $\phi$ and evaluating at $\phi=P_0$, yields the bilinear form
\begin{equation*}
\mathcal{F}_{\phi\phi}(P_{0},\bar{\beta})[x,y]=\left(
\mathcal{H}_{1}[x,y],\ \mathcal{H}_{2}[x,y],\ 0,\ 0,\ \mathcal{H}_{3}[x,y]\right)^{\mathrm{T}},
\end{equation*}
where $y=(y_{1},y_{2},y_{3},y_{4},y_{5})^{\mathrm{T}}$ and
\begin{align*}
\mathcal{H}_{1}[x,y]\ = \ & \frac{\beta_{s}(\theta x_{2}+\epsilon x_{3}+x_{4})[S^{0}(y_{1}+y_{2}+ y_{3}+y_{4}+\int^{\infty}_{0}y_{5}(\tau){\rm d}\tau)-N^{0}y_{1}]}{(N^{0})^2}\\
&\ +\frac{\beta_{s}(\theta y_{2}+\epsilon y_{3}+y_{4})[S^{0}(x_{1}+x_{2}+ x_{3}+x_{4}+\int^{\infty}_{0}x_{5}(\tau){\rm d}\tau)-N^{0}x_{1}]}{(N^{0})^2},\\
\mathcal{H}_{2}[x,y]\ = \ & -\frac{\bar{\beta}(\theta x_{2}+\epsilon x_{3}+x_{4})\int^{\infty}_{0}\beta_{r0}(\tau)[r^{0}(\tau)(y_{1}+y_{2}+ y_{3}+y_{4}+\int^{\infty}_{0}y_{5}(h){\rm d}h)-N^{0}y_{5}(\tau)]{\rm d}\tau}{(N^{0})^2}\\
&\ -\frac{\bar{\beta}(\theta y_{2}+\epsilon y_{3}+y_{4})\int^{\infty}_{0}\beta_{r0}(\tau)[r^{0}(\tau)(x_{1}+x_{2}+ x_{3}+x_{4}+\int^{\infty}_{0}x_{5}(h){\rm d}h)-N^{0}x_{5}(\tau)]{\rm d}\tau}{(N^{0})^2}\\
&\ -\mathcal{H}_{1}[x,y],\\
\mathcal{H}_{3}[x,y]\ = \ & \frac{\bar{\beta}\beta_{r0}(\tau)(\theta x_{2}+\epsilon x_{3}+x_{4})[r^{0}(\tau)(y_{1}+y_{2}+ y_{3}+y_{4}+\int^{\infty}_{0}y_{5}(\tau){\rm d}\tau)-N^{0}y_{5}]}{(N^{0})^2}\\
&\ +\frac{\bar{\beta}\beta_{r0}(\tau)(\theta y_{2}+\epsilon y_{3}+y_{4})[r^{0}(\tau)(x_{1}+x_{2}+ x_{3}+x_{4}+\int^{\infty}_{0}x_{5}(\tau){\rm d}\tau)-N^{0}x_{5}]}{(N^{0})^2}.
\end{align*}
Differentiating $\mathcal F$ with respect to $\phi$ and $\bar\beta$, and evaluating at $\phi=P_0$, yields
\begin{equation*}
\mathcal{F}_{\phi\bar{\beta}}(P_{0},\bar{\beta})x=\left(
0,\ \frac{\int^{\infty}_{0}\beta_{r0}(\tau)r^{0}(\tau){\rm d}\tau(\theta x_{2}+\epsilon x_{3}+x_{4})}{N^{0}},\ 0,\ 0,\ -\frac{\beta_{r0}(\tau)r^{0}(\tau)(\theta x_{2}+\epsilon x_{3}+x_{4})}{N^{0}}\right)^{\mathrm{T}}.
\end{equation*}
Finally, we obtain the following quantities
\begin{align*}
\langle \mathcal{F}_{\phi\phi}(P_{0},\bar{\beta}^{*})[\hat{x},\hat{x}],\hat{\xi}\rangle\ &=\ a,\\
\langle \mathcal{F}_{\phi\bar{\beta}}(P_{0},\bar{\beta}^{*})\hat{x},\hat{\xi}\rangle\ &=\ \frac{1}{N^0}\left(\theta+\frac{(1-\rho)\sigma\epsilon}{\gamma_{A}+u}+\frac{\rho\sigma}{\gamma_{I}+\mu+u}\right)\int^{\infty}_{0}\beta_{r0}(\tau)r^{0}(\tau){\rm d}\tau>0.
\end{align*}
where the expression of $a$ is given in \eqref{EQ_DI4}. Hence, the bifurcation at $(\phi,\bar\beta)=(P_0,\bar\beta^*)$ is transcritical, and the sign of $a$ determines the criticality of the bifurcation, cf., e.g., \cite{kielhofer2012bifurcation,martcheva2020lyapunov}.

\section{Numerical bifurcation analysis}\label{sec3}
\begin{figure}[t!]
\centering
\includegraphics[width=3in,height=2.5in]{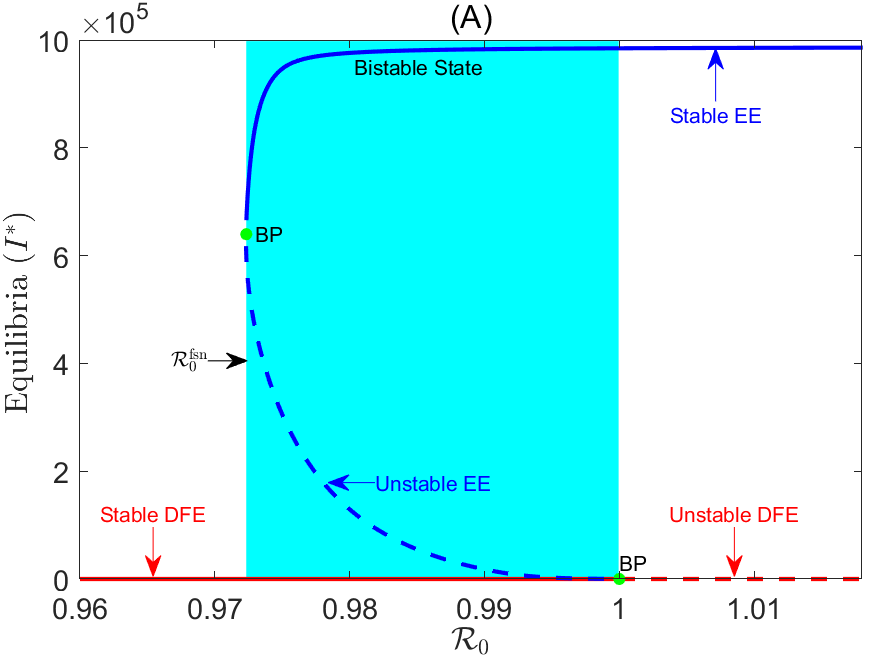}
\hfil
\includegraphics[width=3in,height=2.5in]{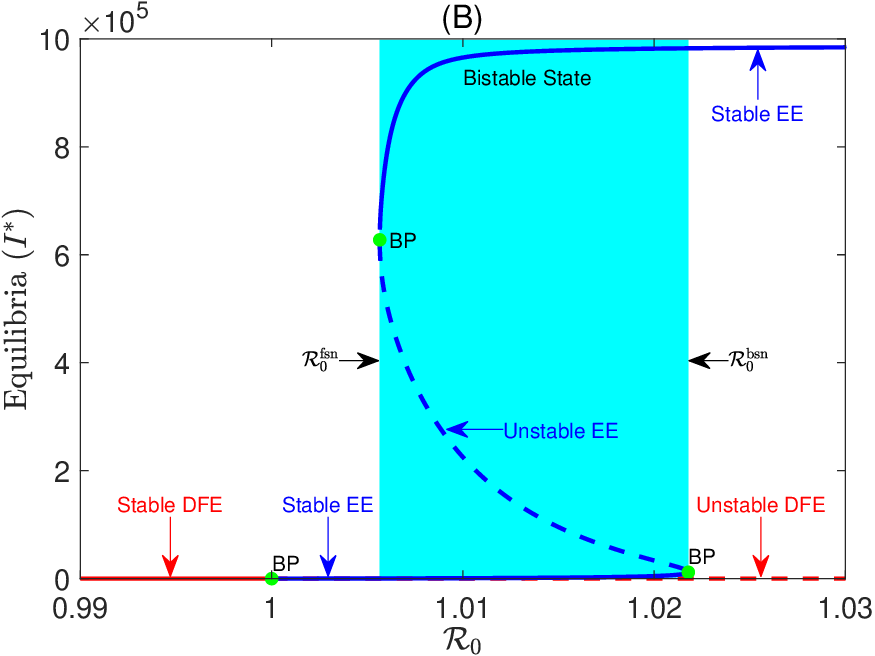}
\caption{
The number of symptomatic infected individuals at equilibria $I(t)=I^*$ varying with $\mathcal{R}_{0}$. (A) $\beta_{s}=0.1$ (i.e., $a=9.7232\times10^{-6}>0$): the DFE undergoes the backward bifurcation at $\mathcal{R}_{0}=1$. (B) $\beta_{s}=0.10345$ (i.e., $a=-1.1558\times10^{-5}<0$): the DFE undergoes the forward bifurcation at $\mathcal{R}_{0}=1$, and the saddle-node bifurcation arises from the stable EE at $\mathcal{R}_{0} =\mathcal R_0^{\rm bsn}\approx1.0218$. Green dots: bifurcation point; cyan regions: the bistable state exists.}
\label{AML_fig1}
\end{figure}

In order to illustrate and corroborate the analytical results, we present some numerical computations. We choose the values of parameters as follows \cite{xue2022infectivity,jing2023vaccine,wu2023prediction}:
$\Lambda=20000$, $\theta=0.55$, $\epsilon=0.55$, $u=1/(75\times365)$, $\alpha=10^{-6}$, $\sigma=1/5.2$, $\rho=0.4$, $\gamma_{A}=1/14$, $\gamma_{I}=1/7$, $\mu=0.02$, and
\begin{equation*}
\beta_{r0}(\tau)=\left\{
\begin{split}
&1-\eta,&&\tau<\hat{\tau},\\
&1-\eta {\rm e}^{-\gamma(\tau-\hat{\tau})},&&\tau\geq \hat{\tau},
\end{split}
\right.
\end{equation*}
where $\hat{\tau}=200$, $\eta=0.2$, and $\gamma=0.5$.
To investigate the impact of $\mathcal{R}_{0}$ on the dynamics of \eqref{COVID_BB_EQ1}, we plot the bifurcation diagram with the horizontal axis $\mathcal R_0$ by varying the bifurcation parameter $\bar\beta$ (due to the continuous dependence of $\mathcal R_0$ on $\bar\beta$), see Figure~\ref{AML_fig1}.

We first consider the case $\beta_{s}=0.1$ (i.e., $a=9.7232\times10^{-6}>0$), the bifurcation from the disease-free equilibrium (DFE) at $\mathcal{R}_{0}=1$ is backward transcritical, and the bifurcating branch extends backward to $\mathcal R_0=\mathcal R_0^{\rm fsn}\approx 0.9723$, at which a saddle-node bifurcation occurs; moreover, for $\mathcal{R}_{0}\in(\mathcal R_0^{\rm fsn},1)$ model \eqref{COVID_BB_EQ1} exhibits the bistable state, in which case a stable endemic equilibrium (EE) and a stable DFE coexist, cf., Figure~\ref{AML_fig1}(A). This indicates, for $\mathcal R_0\in(\mathcal R_0^{\rm fsn},1)$, that the transmission of the infectious disease can be controlled only with sufficiently small initial number of symptomatic infected individuals (SIIs).

Next, we consider the case $\beta_{s}=0.10345$ (i.e., $a=-1.1558\times10^{-5}<0$), the bifurcation from the DFE at $\mathcal{R}_{0}=1$ is forward transcritical, moreover, the bifurcating branch extends forward to $\mathcal{R}_{0}=\mathcal R_0^{\rm bsn}\approx1.0218$, at which the saddle-node bifurcation occurs and its bifurcating branch extends backward to the other saddle-node bifurcation point at $\mathcal R_0=\mathcal R_0^{\rm fsn}\approx 1.0057$, cf., Figure~\ref{AML_fig1}(B). The numerical computations also suggest, that (i) for $\mathcal{R}_{0}\in(1,\mathcal R_0^{\rm fsn})$ the number of SIIs remains small; (ii) for $\mathcal{R}_{0}\in (\mathcal R_0^{\rm fsn},\mathcal R_0^{\rm bsn})$ model \eqref{COVID_BB_EQ1} exhibits the bistable state, in which case two stable EEs coexist, and thus the number of SIIs converges to the lower and upper stable EEs for sufficiently small and large initial data, respectively; (iii) for $\mathcal{R}_{0}>\mathcal R_0^{\rm bsn}$ the number of SIIs always converges to the upper stable EE. These results indicate that $\mathcal{R}_{0}=\mathcal R_0^{\rm fsn}$ is the first critical threshold beyond which the disease can be severe for sufficiently large initial number of SIIs, and $\mathcal{R}_{0}=\mathcal R_0^{\rm bsn}$ is the second critical threshold beyond which the disease becomes severe regardless of the initial number of SIIs.

\section{Conclusion}\label{sec4}
In this paper, we have developed the immune age-structured model \eqref{COVID_BB_EQ1} that takes into account the decline of immunity. Using Lyapunov-Schmidt reduction, we have proved the backward and forward transcritical bifurcations for different values of parameters. We have also presented some numerical computations on various bifurcations, which allow us to explore the nonlinear relations between the immune parameters and the equilibrium of disease. In particular, we found that the saddle-node bifurcation occurs on the extended branch of the forward transcritical bifurcation, and two stable endemic equilibria coexist for some values of parameters.

These results can help us to obtain a deeper understanding of the dynamics of emerging infectious diseases. Incorporating such immunological factors, we are able to gain insights into the conditions that favor disease persistence and severity. We found that the initial number of infected individuals plays a crucial role in determining the disease severity, which emphasizes the importance of early detection and containment. Targeting high-risk individuals and implementing timely interventions may mitigate the impact of emerging infectious diseases and contain their spread.

\section*{Acknowledgements}
L.X. is funded by the National Natural Science Foundation of China 12171116 and Fundamental Research Funds for the Central Universities of China 3072020CFT2402. J.Y. is supported by the Fundamental Research Funds for the Central Universities of China 3072024CFJ2408.

\bibliographystyle{elsarticle-num}
\bibliography{JXY_Dynamics_Immune}

\begin{thebibliography}{10}
\expandafter\ifx\csname url\endcsname\relax
  \def\url#1{\texttt{#1}}\fi
\expandafter\ifx\csname urlprefix\endcsname\relax\def\urlprefix{URL }\fi
\expandafter\ifx\csname href\endcsname\relax
  \def\href#1#2{#2} \def\path#1{#1}\fi

\bibitem{xue2022infectivity}
L.~Xue, S.~Jing, K.~Zhang, R.~Milne, H.~Wang, Infectivity versus fatality of
  {SARS}-{CoV}-2 mutations and influenza, International Journal of Infectious
  Diseases 121 (2022) 195--202.

\bibitem{jing2023vaccine}
S.~Jing, R.~Milne, H.~Wang, L.~Xue, Vaccine hesitancy promotes emergence of new
  {SARS}-{CoV}-2 variants, Journal of Theoretical Biology 570 (2023) 111522.

\bibitem{wangari2024emergence}
I.~M. Wangari, Emergence of a reversed backward bifurcation, reversed
  hysteresis effect, and backward bifurcation phenomenon in a {COVID}-19
  mathematical model, Mathematical Methods in the Applied Sciences 47~(4)
  (2024) 2250--2272.

\bibitem{martcheva2015introduction}
M.~Martcheva, An Introduction to Mathematical Epidemiology, Vol.~61, Springer,
  New York, 2015.

\bibitem{martcheva2020lyapunov}
M.~Martcheva, H.~Inaba, A {L}yapunov--{S}chmidt method for detecting backward
  bifurcation in age-structured population models, Journal of Biological
  Dynamics 14~(1) (2020) 543--565.

\bibitem{castillo2004dynamical}
C.~Castillo-Chavez, B.~Song, Dynamical models of tuberculosis and their
  applications, Mathematical Biosciences and Engineering 1~(2) (2004) 361--404.

\bibitem{yang2022backward}
J.~Yang, M.~Zhou, X.~Li, Backward bifurcation of an age-structured epidemic
  model with partial immunity: the {L}yapunov--{S}chmidt approach, Applied
  Mathematics Letters 133 (2022) 108292.

\bibitem{rahman2022covid}
S.~Rahman, M.~M. Rahman, M.~Miah, M.~N. Begum, M.~Sarmin, M.~Mahfuz, M.~E.
  Hossain, M.~Z. Rahman, M.~J. Chisti, T.~Ahmed, et~al., {COVID}-19
  reinfections among naturally infected and vaccinated individuals, Scientific
  Reports 12~(1) (2022) 1438.

\bibitem{risk2022covid}
M.~Risk, S.~S. Hayek, E.~Schiopu, L.~Yuan, C.~Shen, X.~Shi, G.~Freed, L.~Zhao,
  {COVID}-19 vaccine effectiveness against omicron ({B}.1.1.529) variant
  infection and hospitalisation in patients taking immunosuppressive
  medications: a retrospective cohort study, The Lancet Rheumatology 4~(11)
  (2022) e775--e784.

\bibitem{van2002reproduction}
P.~van~den Driessche, J.~Watmough, Reproduction numbers and sub-threshold
  endemic equilibria for compartmental models of disease transmission,
  Mathematical Biosciences 180~(1-2) (2002) 29--48.

\bibitem{diekmann1990definition}
O.~Diekmann, J.~A.~P. Heesterbeek, J.~A. Metz, On the definition and the
  computation of the basic reproduction ratio ${R}_{0}$ in models for
  infectious diseases in heterogeneous populations, Journal of {M}athematical
  {B}iology 28~(4) (1990) 365--382.

\bibitem{kielhofer2012bifurcation}
H.~Kielh{\"o}fer, Bifurcation Theory: An Introduction with Applications to
  Partial Differential Equations, Springer, New York, 2012.

\bibitem{wu2023prediction}
Y.~Wu, W.~Zhou, S.~Tang, R.~A. Cheke, X.~Wang, Prediction of the next major
  outbreak of {COVID}-19 in {M}ainland {C}hina and a vaccination strategy for
  it, Royal Society Open Science 10~(8) (2023) 230655.

\end{thebibliography}
\end{document}